\documentclass{ifacconf}
\usepackage[utf8]{inputenc}

\usepackage{graphicx}   
\usepackage{epsfig} 
\usepackage[T1]{fontenc} 
\usepackage{amssymb} 
\usepackage{amsmath} 
\usepackage{algorithm}
\usepackage[noend]{algpseudocode}
\usepackage{graphicx}
\usepackage{lipsum}%
\usepackage{xcolor}
\usepackage{caption}
\usepackage{float}
\usepackage{enumitem}
\makeatletter
\let\old@ssect\@ssect 
\makeatother

\usepackage{natbib}
\usepackage{hyperref}
\makeatletter
\def\@ssect#1#2#3#4#5#6{%
  \NR@gettitle{#6}
  \old@ssect{#1}{#2}{#3}{#4}{#5}{#6}
}
\makeatother

\begin{document}
\begin{frontmatter}
\title{SECREDAS: Safe and (Cyber-)Secure Cooperative and Automated Mobility}
\author[First]{Chris van der Ploeg}
\author[First]{Jacco van de Sluis}
\author[Second]{Sebastian Gerres}
\author[Third]{Szabolcs Novaczki}
\author[Third]{Andr{\'a}s Wippelhauser}
\author[Fourth]{Eric Nassor}
\author[Fourth]{Julien Sevin}
\author[Fifth]{Andr{\'a}s Gazdag}
\author[Fifth]{Gergely Bicz{\'o}k}

\address[First]{TNO Integrated Vehicle Safety, Automotive Campus 30, 5708 JZ, Helmond, The Netherlands}
\address[Second]{Merantix Momentum GmbH, Max-Urich-Straße 3, 13355, Berlin, Germany}
\address[Third]{Commsignia, Irinyi Jozsef u. 4-20, 1117, Budapest, Hungary}
\address[Fourth]{Canon Research Centre France, Rue de la Touche Lambert, 35510, Cesson-Sevigne, France}
\address[Fifth]{CrySyS Lab, Dept. of Networked Systems and Services, Budapest Univesrity of Technology and Economics, Műegyetem rkp. 3, 1111, Budapest, Hungary}

\begin{abstract}
    Infrastructure-to-Vehicle (I2V) and Vehicle-to-Infrastructure (V2I) communication is likely to be a key-enabling technology for automated driving in the future. Using externally placed sensors, the digital infrastructure can support the vehicle in perceiving surroundings that would otherwise be difficult to perceive due to, for example, high traffic density or bad weather. 
    Conversely, by communicating on-board vehicle measurements, the environment can more accurately be perceived in locations which are not (sufficiently) covered by digital infrastructure. The security of such communication channels is an important topic, since malicious information on these channels could potentially lead to a reduction in overall safety. Collective perception contributes to raising awareness levels and an improved traffic safety. 
    In this work, a demonstrator is introduced, where a variety of novel techniques have been deployed to showcase an overall architecture for improving vehicle and vulnerable road user safety in a connected environment. The developed concepts have been deployed at the Automotive Campus intersection in Helmond (NL), in a field testing setting.
\end{abstract}
\end{frontmatter}
\begin{figure*}[t]
    \centering
    \includegraphics[width=0.9\textwidth]{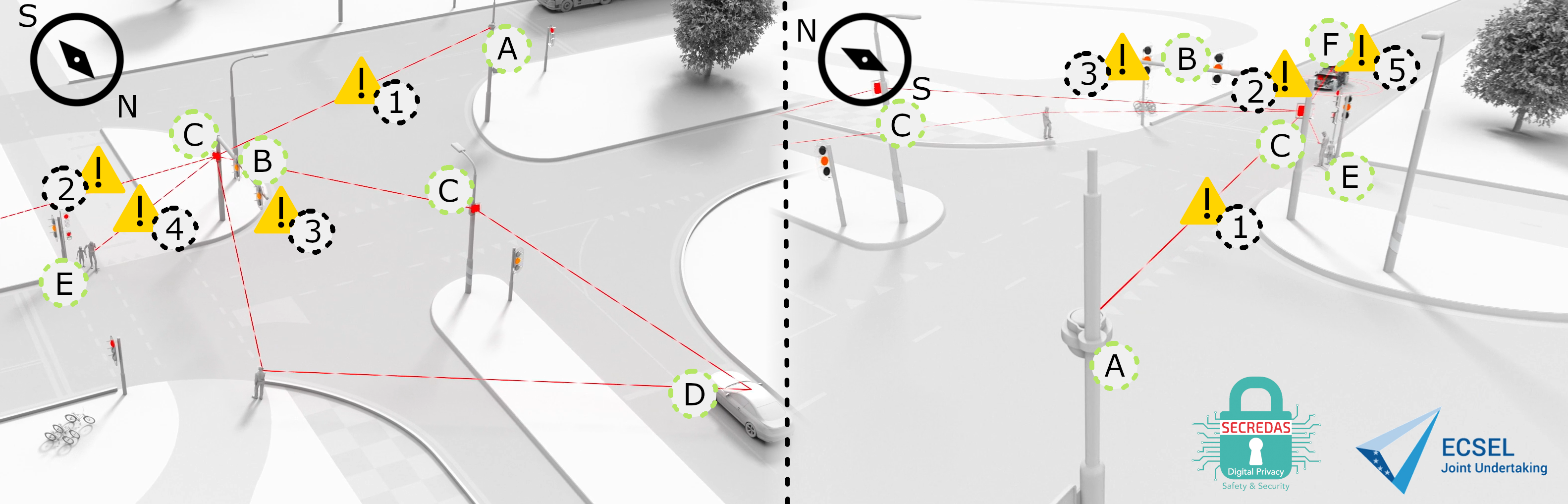}
    \caption{Visual representation of all considered cyber-attacks and actors.}
    \label{fig:secredas_attacks}
\end{figure*}
\section{Introduction}
{Cooperative Connected and Automated Mobility (CCAM) is an overall term for smart mobility solutions enabled by communication technology, to improve traffic safety and to achieve more efficient, comfortable and sustainable ways of transport. Example applications of connected vehicles (e.g., the Grand Cooperative driving Challenge~\cite{7553038}) have shown the benefits of CCAM to the wide public. A key-enabling technology for CCAM is the digital infrastructure for automated driving functions~\cite{Silva2017} to raise the level of awareness by sharing perception information in mixed traffic conditions. Using this technology, perception is brought to a collective problem rather than a single-vehicle problem, ensuring that a traffic participant is able to anticipate further ahead and hence increase single-vehicle and overall traffic safety~\cite{Shan20211}. Although these technologies can bring great benefits, the use of a connected platform introduces inherent vulnerabilities (e.g., proneness to cyber-attacks~\cite{Petit2015546}). The SECREDAS project (stands for "Product Security for Cross Domain Reliable Dependable Automated Systems"~\cite{SECREDAS}) focused on cyber-security solutions for CCAM, to develop use cases, real-life threat scenarios, a security framework and architecture, and to use common technology elements to be integrated into demonstrator scenarios. The work presented in this paper focuses on securing a cooperative road intersection use case, with automated vehicles and infrastructure elements sharing perception information to improve situational awareness. Vehicle-to-infrastructure (V2I, I2V) communication is used to exchange sensor information amongst vehicles and digital infrastructure. Multiple attack scenarios are identified and used to developed detection and mitigation solutions. The five attack scenarios to be executed are:
\begin{enumerate}
	\item \textbf{Malicious Cooperative Perception Messages (CPM):} The roadside sensor's communication is hacked broadcasting malicious CPMs, containing false object data about detected traffic;
	\item \textbf{Hacked vehicle communication:} A hijacked automated vehicle broadcasting malicious cooperative awareness messages (CAMs) indicating that it is stopping for a red light, while in reality it will violate the red light;
	\item \textbf{Hacked traffic light communication:} The traffic lights are hacked and broadcasting malicious Signal Phase and Timing (SPaT) messages indicating that all lights are green;
	\item \textbf{Protection of vulnerable road users (VRU):} A pedestrian carrying an ultra wide-band (UWB) tag and a smartphone is crossing the intersection with an approaching (hijacked) vehicle;
	\item \textbf{Hacked on-board camera for object detection:} A camera sensor of the automated vehicle is hacked and an adversarial image attack is performed to impair the vehicle's object detection.
\end{enumerate}
The hardware and software, needed to safely and securely detect and mitigate these attacks, is fully integrated and deployed at the Automotive Campus intersection in Helmond in field testing for validation of the use cases, proof-of-concept of the threat mitigation solutions and demonstration in real-life environment. 
The demonstrator use case is an intersection crossing scenario consisting of smart digital infrastructure (sensors, I2V communication), automated vehicles (sensors, V2I communication) and vulnerable road users (VRUs).
Fig.~\ref{fig:secredas_attacks} depicts the main actors and elements for the intersection use case and the used attack scenarios. The cooperative intersection is equipped with a roadside camera (A) traffic lights and control (B) and road side units (RSUs) for I2V communications (C), all part of a traffic management system. Other actors in this scenario are the victim vehicle (D) which may suffer consequences of a considered attack, the VRUs (E) and a connected automated vehicle (F), both communicating with an on-board unit (OBU). This setting is the starting point of the use case supporting multiple attack scenarios.
The identified attack scenarios are used to further explain the developed detection and mitigation solutions (in the next section). All experimental tests have been designed for, and done with, a driver-in-the-loop, i.e., the driver is responsible for maintaining the safety at all times. Any benefiting technologies that arise from this project are therefore proposed as driver-assistance features.}
\section{Overview and Challenges}
 In vehicle-to-anything (V2X) communication, three security principles are observed: authentication, integrity protection and non-repudiation. It is imperative that we can identify the sender of an erroneous/fake message, therefore, the sender should not be able to deny sending it; otherwise, we could not exclude the non-compliant participants from the system. To this end, the current standard defines a Public Key Infrastructure (PKI), which is similar to the certificate-based system used in the Internet. Every vehicle receives an asymmetric key-pair and a certificate signed by a trusted entity; after this, the vehicle signs all its messages with its private key and sends them together with its certificate. At the other end, the receiving party can check the sender's identity and the integrity of the message using the public key in the certificate. Combined with the use of short-lived pseudonyms to ensure location privacy, and checking the chain of Certificate Authorities (CA) up until the Root CA, the mechanism is robust and ensures authentication, message integrity and non-repudiation. The mechanism makes it difficult to send out self-constructed messages, as the keys, needed for a valid certificate and signature, are stored in a special, tamper-resistant unit, the Hardware Security Module (HSM).
 
Although it is reasonably difficult to get a valid signature on a self-fabricated message, sending falsified messages is far from impossible. If there is a vulnerability in the on-board software system of a vehicle, which can be exploited to take control of the data flowing through the Controller Area Network (CAN) bus, then the adversary can overwrite the exact values the OBU receives from the vehicle sensors. Not being aware of the attack, the OBU happily puts the valid signature on the falsified values, before sending out the seemingly intact, but content-wise modified message towards the world. It can be seen that V2X security extends beyond the control of V2X stack providers: however secure their own system is, carmakers/integrators also bear significant responsibility in securing V2X communications. Anyway, 100\% software security is impractical, and vulnerabilities are bound to be found. 

The variety of misbehaviours/attacks that can be instantiated by modifying information received by an OBU (and sent out in V2X messages) is enormous. There are simpler ones, e.g., where an attacker can forge Decentralized Environmental Notification Messages (DENM)  indicating a traffic jam to clear the path for herself. Also, an attacker can swiftly modify the geo-coordinates in a CAM message in a way that her car appears virtually very close to another vehicle, triggering its emergency braking mechanism. While attacks like these can endanger vehicles and drivers alike, they are relatively straightforward. Unfortunately, there are more intricate possibilities. During a sybil attack, the adversary sends fake CAM messages from one or more ghost vehicles. By emulating multiple sybil cars, whole fleets can be faked, which can even affect traffic management mechanisms. Note that unintended failures may also happen at the sensor-level, e.g., an erroneous GPS module may produce inconsistent readings, or a malfunctioning accelerometer could trigger the emergency brakes. Furthermore, sensors can be tricked maliciously, e.g., the close proximity of a bag of ice may make our vehicle send out DENM messages indicating slippery road surface in the middle of the summer.

It is clear that the detection of such misbehavior can be done in a variety of ways. Here, we advocate for a quick, simple and extensible first-line, in-car mechanism with four key characteristics. First, it should work autonomously without trusting another entity. Second, it should be based on quick and effective filtering methods as opposed to computation-heavy, complex algorithms to enable a quick reaction. Third, confidence intervals and dynamic thresholds should be used to improve detection accuracy and minimize false positives. Last, as V2X communication itself is a developing technology, the misbehavior detection mechanism should be flexible to accommodate new standards and incorporate future algorithmic improvements. 

In the V2X systems the application information shared by a V2X node is limited by the service specific permission field which is part of the certificates. Apart from the vehicle on-board sensors the road side infrastructure can also be prone to cyber-physical attacks. The roadside module is allowed to advertise regulatory road signs or traffic light information. In the case of misbehaviours the false speed limits, road works, traffic directions and lane topologies can be advertised. Maybe the most dangerous scenario that can occur is the modification of traffic light sequences. Practically this can be achieved via malicious SPaT messages. A misbehaviour detection method can facilitate the lane topology and the signal phase sequences in order to evaluate that the traffic lights do not advertise conflicting movement allowances.
\section{Proposed Architecture}
\begin{figure}[t]
    \def\svgwidth{1\columnwidth}
\begingroup%
  \makeatletter%
  \providecommand\color[2][]{%
    \errmessage{(Inkscape) Color is used for the text in Inkscape, but the package 'color.sty' is not loaded}%
    \renewcommand\color[2][]{}%
  }%
  \providecommand\transparent[1]{%
    \errmessage{(Inkscape) Transparency is used (non-zero) for the text in Inkscape, but the package 'transparent.sty' is not loaded}%
    \renewcommand\transparent[1]{}%
  }%
  \providecommand\rotatebox[2]{#2}%
  \newcommand*\fsize{\dimexpr\f@size pt\relax}%
  \newcommand*\lineheight[1]{\fontsize{\fsize}{#1\fsize}\selectfont}%
  \ifx\svgwidth\undefined%
    \setlength{\unitlength}{442.20472441bp}%
    \ifx\svgscale\undefined%
      \relax%
    \else%
      \setlength{\unitlength}{\unitlength * \real{\svgscale}}%
    \fi%
  \else%
    \setlength{\unitlength}{\svgwidth}%
  \fi%
  \global\let\svgwidth\undefined%
  \global\let\svgscale\undefined%
  \makeatother%
  \begin{picture}(1,1.23076923)%
    \lineheight{1}%
    \setlength\tabcolsep{0pt}%
    \put(0,0){\includegraphics[width=\unitlength,page=1]{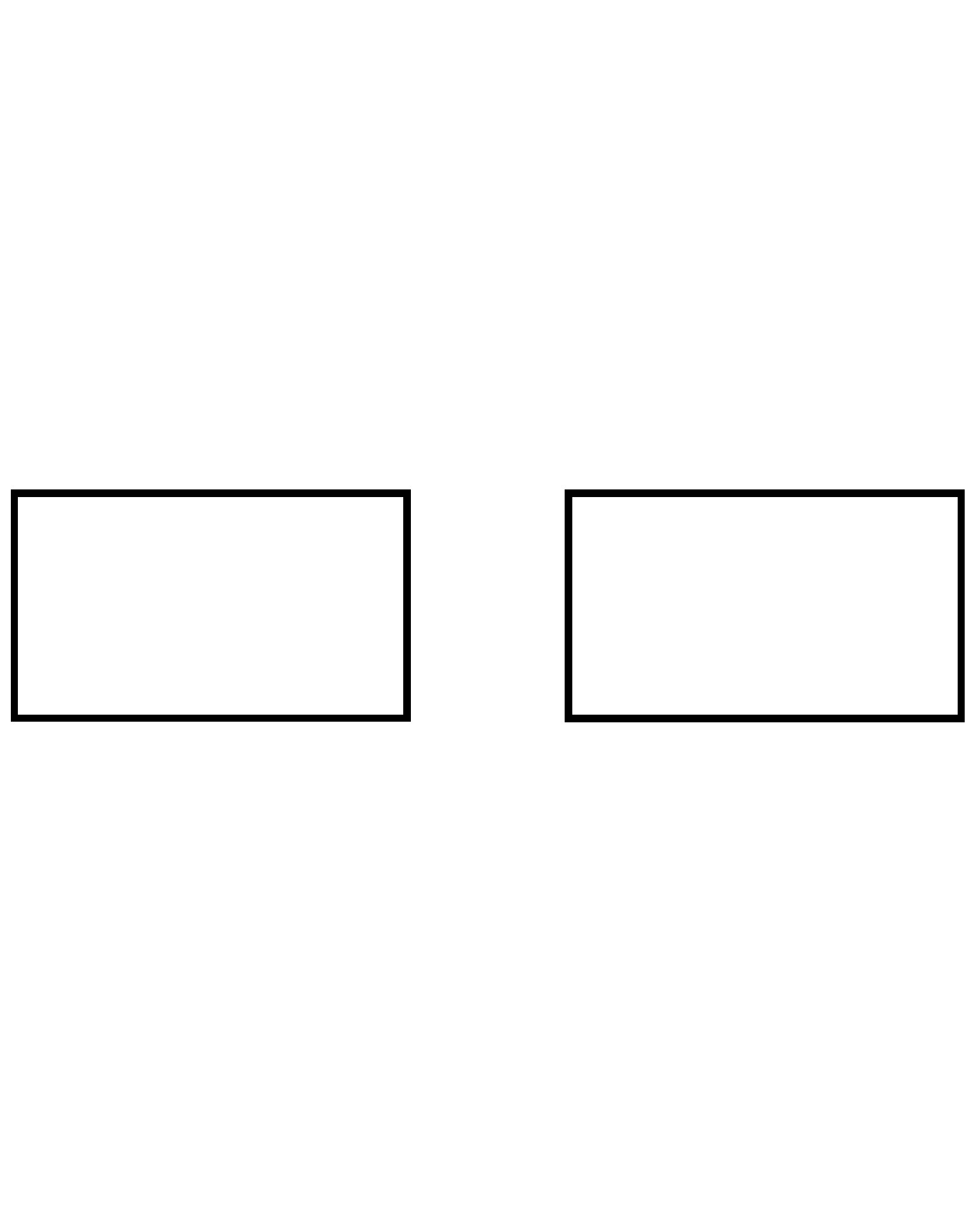}}%
    \put(0.78439962,0.62133339){\color[rgb]{0,0,0}\makebox(0,0)[t]{\lineheight{1.25}\smash{\begin{tabular}[t]{c}Traffic light\\control\end{tabular}}}}%
    \put(0.06542788,0.59481704){\color[rgb]{0,0,0}\makebox(0,0)[lt]{\lineheight{1.25}\smash{\begin{tabular}[t]{l}Roadside unit\end{tabular}}}}%
    \put(0,0){\includegraphics[width=\unitlength,page=2]{fig3.pdf}}%
    \put(0.78354551,0.14113382){\color[rgb]{0,0,0}\makebox(0,0)[t]{\lineheight{1.25}\smash{\begin{tabular}[t]{c}Connected\\VRU\end{tabular}}}}%
    \put(0,0){\includegraphics[width=\unitlength,page=3]{fig4.pdf}}%
    \put(0.22399487,1.10615315){\color[rgb]{0,0,0}\makebox(0,0)[t]{\lineheight{1.25}\smash{\begin{tabular}[t]{c}Roadside\\camera\end{tabular}}}}%
    \put(0.21944587,0.13843918){\color[rgb]{0,0,0}\makebox(0,0)[t]{\lineheight{1.25}\smash{\begin{tabular}[t]{c}Connected \\vehicle OBU\end{tabular}}}}%
    \put(0,0){\includegraphics[width=\unitlength,page=4]{fig5.pdf}}%
    \put(0.37502,0.35732299){\color[rgb]{0,0,0}\makebox(0,0)[lt]{\lineheight{1.25}\smash{\begin{tabular}[t]{l}CAM\end{tabular}}}}%
    \put(0.2255685,0.83382372){\color[rgb]{0,0,0}\makebox(0,0)[lt]{\lineheight{1.25}\smash{\begin{tabular}[t]{l}CPM\end{tabular}}}}%
    \put(0.08796204,0.83382377){\color[rgb]{0,0,0}\makebox(0,0)[lt]{\lineheight{1.25}\smash{\begin{tabular}[t]{l}DENM\end{tabular}}}}%
    \put(0.44215687,0.62239671){\color[rgb]{0,0,0}\makebox(0,0)[lt]{\lineheight{1.25}\smash{\begin{tabular}[t]{l}SPaT/\\MAP\end{tabular}}}}%
    \put(0.42499738,0.7546017){\color[rgb]{0,0,0}\makebox(0,0)[lt]{\lineheight{1.25}\smash{\begin{tabular}[t]{l}Light\\control\end{tabular}}}}%
    \put(0.22477637,0.36076531){\color[rgb]{0,0,0}\makebox(0,0)[lt]{\lineheight{1.25}\smash{\begin{tabular}[t]{l}CPM\end{tabular}}}}%
    \put(0.08899186,0.35910748){\color[rgb]{0,0,0}\makebox(0,0)[lt]{\lineheight{1.25}\smash{\begin{tabular}[t]{l}DENM\end{tabular}}}}%
    \put(0.37360167,0.83267647){\color[rgb]{0,0,0}\makebox(0,0)[lt]{\lineheight{1.25}\smash{\begin{tabular}[t]{l}CAM\end{tabular}}}}%
    \put(0,0){\includegraphics[width=\unitlength,page=5]{fig6.pdf}}%
    \put(0.55425551,0.35124833){\color[rgb]{0,0,0}\makebox(0,0)[lt]{\lineheight{1.25}\smash{\begin{tabular}[t]{l}UWB\end{tabular}}}}%
    \put(0.79608628,0.41965441){\color[rgb]{0,0,0}\makebox(0,0)[t]{\lineheight{1.25}\smash{\begin{tabular}[t]{c}Security\\notification\end{tabular}}}}%
  \end{picture}%
\endgroup%

    \caption{Dataflow architecture of the SECREDAS demonstration}
    \label{fig:architecture2}
\end{figure}
In the scope of the SECREDAS project we built up a software and hardware infrastructure in order to test and evaluate the various misbehaviour algorithms and scenarios. The complete architecture contains the following subsystems (as depicted in Fig.~\ref{fig:architecture2}).

\begin{enumerate}
    \item The road infrastructure subsystem including a controllable traffic light and an RSU.
    \item A roadside perception subsystem, perceiving objects which are traveling within the field of view of the camera.
    \item The vehicle subsystem including an OBU and on-board object detection.
    \item Handheld devices, carried by vulnerable road users along the infrastructure.
\end{enumerate}

The infrastructure perception services subsystem contains a camera system covering the whole intersection. An artificial intelligence based algorithm detects the objects on the camera image and feeds the object position information to the RSU in the form of a CPM. The intersection is also equipped with an UWB subsystem which localizes the pedestrians with suitable smart handheld devices and could warn these pedestrians in case of danger of collision with other actors. The RSU collects the CPM messages and distributes them along the connected actors. Using the CPM generation method we are able to integrate the various sensor sources into the existing V2X architecture while respecting existing security concepts.
In the road infrastructure,  consisting of a RSU and a traffic light, the RSU is able to share the status of the traffic light (i.e., SPaT/MAP messages) to connected road users. Moreover, it acts as a gateway between the V2X networks and the perception services for intercommuncation of DENM, CPM and CAM messages.
The vehicle subsystem contains an OBU and object perception sensors (e.g., radar, camera). The OBU performs the regular V2X communication methods, like CA services, and also deals with the perception sharing services, i.e., the communication of CAM, CPM and DENM messages, which effectively links the object database of the RSU and the OBU.
Vulnerable road users, e.g., pedestrians are vulnerable at a connected intersection since certain misbehaviors of connected and automated vehicles may impose hazards. In the demonstration we assume that the present VRUs carry handheld devices with UWB capabilities, which allows the RSU to locate the pedestrians on the intersection. A return channel is used to warn the VRU if a collision with another road user (connected or not) is imminent.

In the scope of this project the possible malicious attacks against cooperative and automated driving were examined. We consider the misbehaviour detection algorithms to be deployed in accordance with the Cooperative Intelligent Transport Systems and Services (C-ITS) architecture standards. This means that the misbehaviour functionality is a cross-layer entity with presence in all layers of communication. Misbehaviour detection algorithms are deployed next to each functional element in the system. The misbehaviour detection algorithms are connected to an API which is designed to mitigate the misbehaviour events. This API's, called the 'Security Notification API', main goal is to provide a centralized functionality, which limits the error propagation within the V2X network. In the case of a detected malicious attack on the automated vehicle, the Security Notification API is able to revoke the trust from the V2X communication unit by purging its own private keys. In the case of misbehaviour detection of a remote V2X station, the Security Notification API is able to increase the awareness of other vehicles by sending appropriate notification messages to be displayed on the vehicle HMI.
\section{Demonstrators: attacks and subsequent mitigation strategies}
In this section we further elaborate on each of the attack scenarios, the hazard imposed and the detection of such a hazard and, subsequently, the mitigation.
\subsection{Malicious CPM messages}~\label{sec:scen1}
Object detection (e.g., for motion planning or adaptive cruise control) is an example where infrastructural sensor information could be used to the benefit of an automated vehicle. A roadside unit, part of the infrastructure, may be able to perceive objects that are outside the field-of-view of internal vehicle sensors. This information can be used to better anticipate in certain safety critical situations (e.g., a cyclist approaching an intersection while the vehicle is approaching the same intersection at high velocity). Despite its advantages, using this information in closed-loop could result in safety hazards, for example, when the CPM contains malicious information about an object, whereas other on-board sensors state contradicting information. It is up to the vehicle to compile this information and conclude whether a potential security-breach has occurred. In the work, proposed in~\cite{TNOECC21}, an anomaly detection algorithm has been proposed (of which the architecture is depicted in~\cite[Fig. 5]{TNOECC21}, and fully elaborated in~\cite{TNOECC21}), which detects anomalies of different signatures from different sources in a model-driven approach using an extended Kalman filter, describing the dynamics of the VRU. Once a cyber-security breach has been detected (i.e., malicious velocity data), the vehicle sends out a DENM message to the RSU. Subsequently, this V2X communication warning (DENM message) is communicated to other vehicles. As a result, the warning is shown on the driver’s HMI inside a second vehicle which is also approaching the intersection, and on the HMI of the RSU. With this mitigation action in place, other road users are alerted and potentially unsafe reactions to malicious object data are avoided.
\subsection{Hacked vehicle}~\label{sec:scen2}
This scenario assumes the threat of a hacker, fully taking over control of an automated vehicle. The vehicle is assumed to be out of control of any of the standard mitigations systems. The attacker forces the vehicle to drive through a red traffic light of an intersection, thus endangering the vehicle occupants as well as other nearby traffic participants. The attacker also manipulates the communication of the vehicle, such that it broadcasts CAM messages with the intent to stop (deceleration and eventually speed of zero) even though it is not. If this were not the case, the traffic and the road side units at the intersection would be able to anticipate a red-light violation. The faulty CAM messages ensure that also external systems are misled by the hacker.
In order to detect anomalies related to hacked internal vehicle signals, five dedicated types of anomaly detection algorithms have been considered, each one being successively launched at the reception of an ITS message: 
\begin{itemize}
    \item A basic security level anomaly detection module based on local security checks as specified in the standards~\cite{ETSI103097} and~\cite{ETSI1030962}. It corresponds to the first level of security.
    \item A local PHY sensors anomaly detection module based on measurements captured by local PHY sensors of the ego ITS station. It can be applied for both communication technologies considered by ETSI, ITS-G5 (based on IEEE 802.11p) and cellular technologies (4G/5G).
    \item An implausibility anomaly detection based on the payload implausibility of the received ITS messages. Some plausibility detectors are given by standard [ETSI TR 103 460] for CAM messages.
    \item An inconsistency anomaly detection based on the payload implausibility of successive received ITS messages. This is the level 2 of the anomaly detection described in standard~\cite{ETSI103460} for CAM messages.
    \item An enhanced inconsistency anomaly detection module based on the inconsistency between the received message and the information retrieved from the on-board sensors (level 4 of standard~\cite{ETSI103460}) or the information from previous received ITS messages (level 3 of~\cite{ETSI103460}).
\end{itemize}
In the case that an anomaly is detected by one of these anomaly detection modules, an anomaly report is generated (in the form of a DENM for instance), indicating the type of anomaly and the corresponding evidence. Otherwise, the ITS message is considered safe and the information reported by it is injected in the local dynamic map of the ego station.

A first implementation of this architecture has been developed and tested. In particular, a specific anomaly detection algorithm has been specified in order to detect inconsistency anomalies from the reception of CPMs by exploiting the additional redundancy information provided by these messages.  
The developed anomaly detection algorithm allows to detect three types of anomalies involving CPMs.  A position usurpation when two different ITS stations advertise by a CAM or a CPM the same position with two identities for which the on-board sensors have detected the presence of an object, a ghost anomaly when an ITS station advertises the position of a station by a CPM for which the local on-board sensors have detected no vehicle and a hijacked vehicle anomaly when an ITS station advertises no ITS station at a position for which the local on-board sensors have detected a vehicle.
\begin{figure}[t]
    \centering
    \includegraphics[width=\columnwidth]{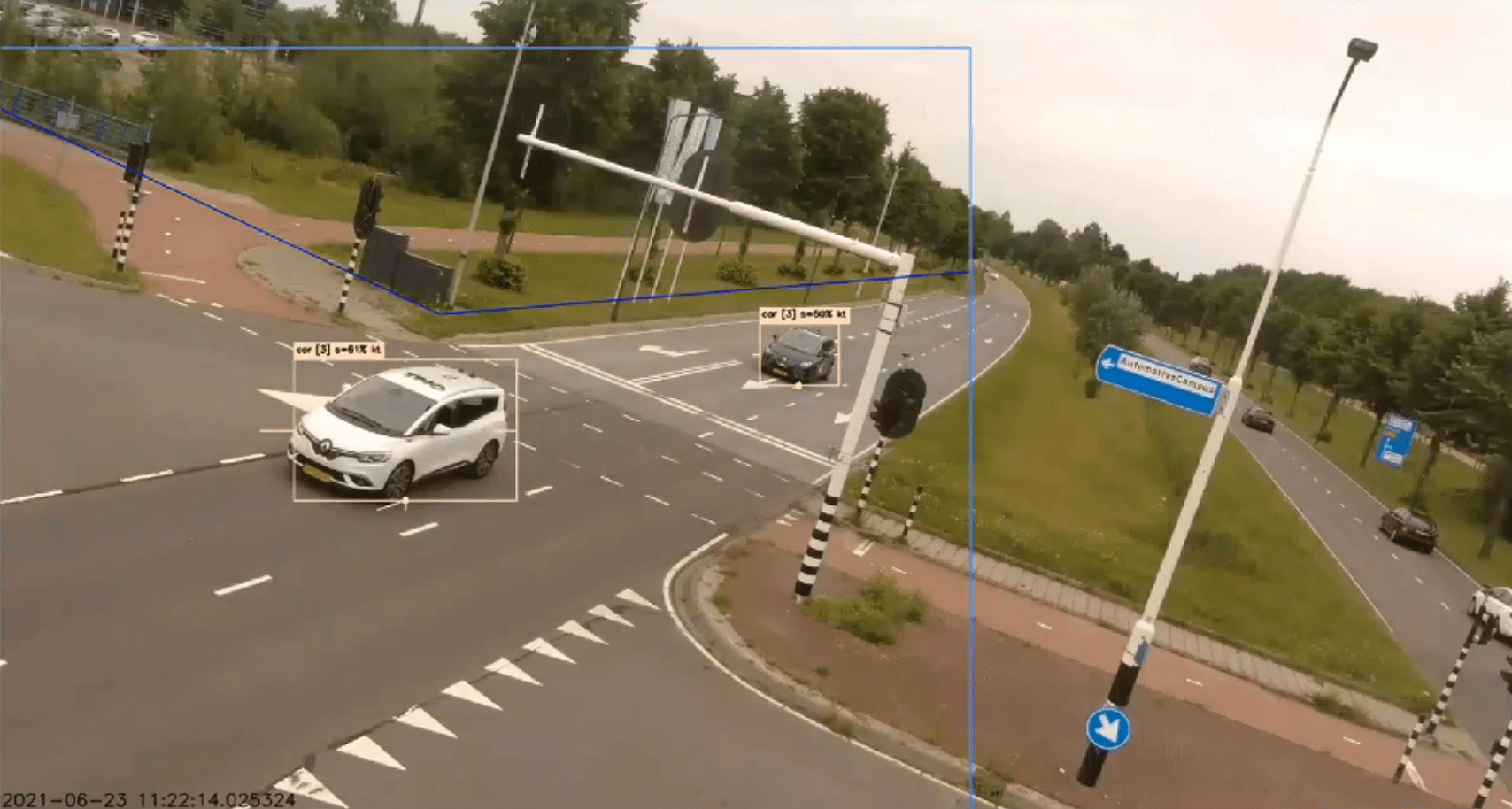}
    \caption{Video feed of the test vehicles observed by the roadside camera.}
    \label{fig:CRF_2}
\end{figure} 
Such anomalies are quite new and in particular, they have not yet been considered in the standard relative to misbehavior reporting service~\cite{ETSI103759} which is currently specified. At the current stage, only anomalies based on CAM receptions have been taken into account and the next step should be the management of the other ITS messages as DENM, VAM and CPM. This additional awareness coupled with the specified global architecture should enhance the anomalies detection process and allows it to integrate definitely the existing ITS security architecture specified in standard~\cite{ETSI102731} which provides security and verification required to secure V2X communications. Once a malicious CAM message is detected, a warning in the form of DENM can be received by other vehicles and is subsequently shown on the HMI in the victim vehicle as it is approaching the intersection. 
A second implemented mitigating action in this scenario is that the RSU requests the traffic light management system to turn the traffic lights to red/orange blinking. This way the occupant of the victim vehicle also sees traffic light phase change from green to red.
\subsection{Hacked traffic lights}~\label{sec:scen3}
This scenario assumes that a hacker is attacking the road infrastructure, in this case the traffic light management system. By gaining control over the traffic lights’ phase status (red, yellow, green) or - at least - over the phase information broadcast by the traffic light communication system, the considered type of attack is that the intersection is broadcasting that all lights are showing a green signal although, in reality, they could be red or orange. Approaching traffic from all sides would assume that it is safe to cross due this green light signal. However, in reality, none of the sides are safe and appropriate mitigating action must be undertaken to maintain safety. The threat is detected by a traffic light anomaly detection algorithm. It takes as input the SPaT/MAP messages sent out by the traffic light system on the intersection. It performs validity checks on the broadcast information and so can detect if the information has been tampered with. In our implementation, the algorithm determines that contradicting traffic lights show a green signal (e.g. all traffic lights on the intersection are suddenly green). The RSU then sends out a warning in the form of a DENM message to all nearby road users. As a second mitigation action, the RSU attempts a request for all traffic lights to turn their signal to red/orange blinking. In reality, the latter cannot be guaranteed since the traffic light system itself is (also) under attack – and the possibilities depend on the type of attack. With these two mitigating actions in place, possible collisions between the intersection users are avoided.
\begin{figure}[t]
    \centering
    \includegraphics[width=0.8\columnwidth]{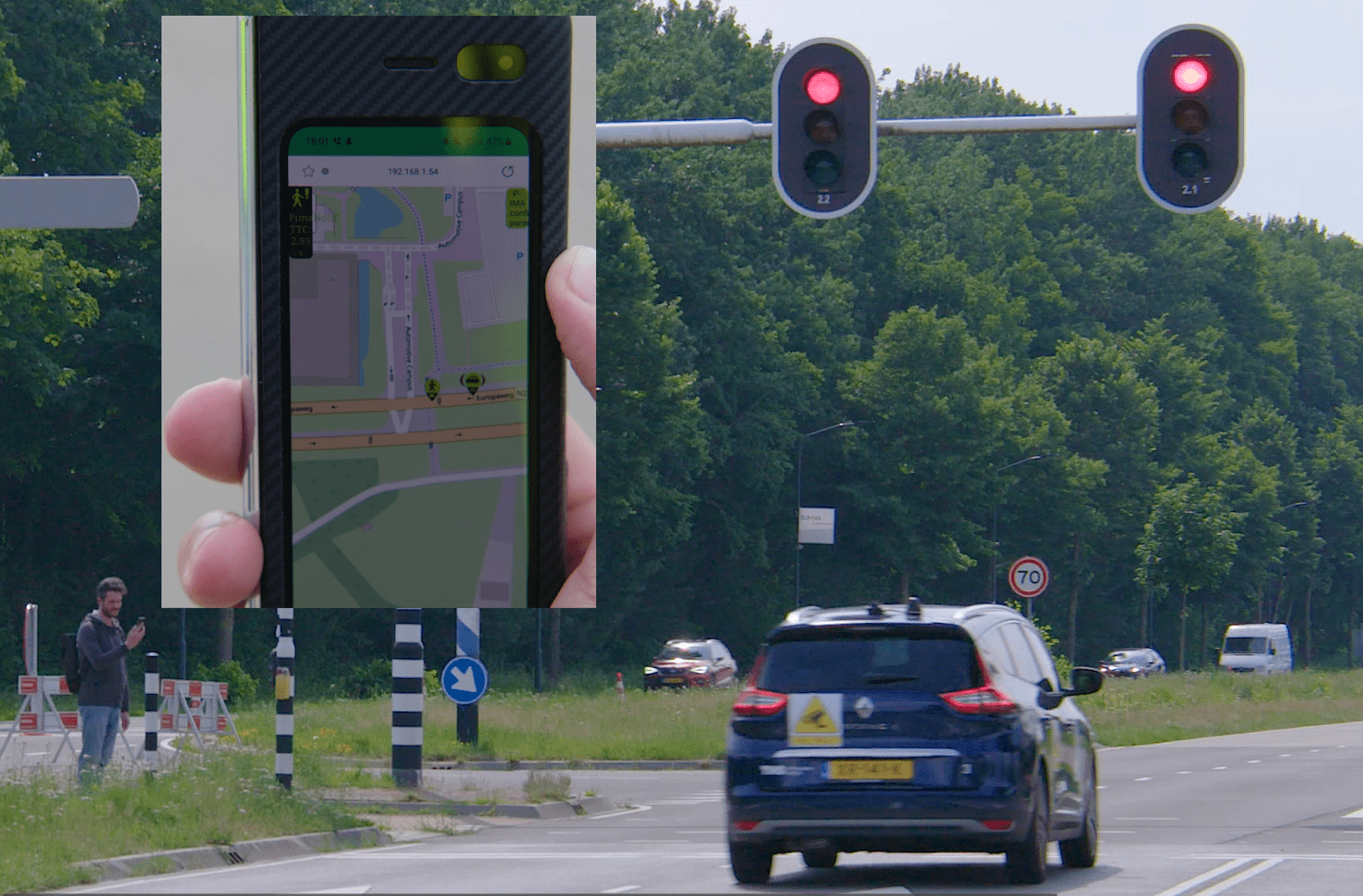}
    \caption{An alerted pedestrian waiting for the hacked vehicle, passing through a red light.}
    \label{fig:CRF_2}
\end{figure} 
\subsection{Protection of pedestrians}~\label{sec:scen4}
The VRU protection is also a key aspect of road safety, since the pedestrians or cyclists do not benefit from the protection provided by the passive safety features of the vehicles. Thus the impact of a possible collision can be more severe for them. The redundancy provided by the collective perception service shared via V2X communication can largely improve the safety of the roads. The smart device of the pedestrian is able to connect to the RSU which detects the approaching vehicle via the roadside camera or via the CPM or CAM messages of neighboring V2X equipped devices. The RSU is also able to identify the approaching malicious vehicle by detecting the misbehaviour using the sensor perception information and the CAMs sent by the vehicle. For example, if the vehicle simulates a gentle stop before the red light in its CAMs, but in fact it proceeds and crosses the red light, the RSU is able to detect the malicious behaviour and notify the VRUs nearby. The malicious behaviour can be detected using the deviation between the CAM and the perception data. In this use case the precise positioning of the vulnerable road users is essential to avoid false positive notifications due to their unpredictable behaviour and small size. This precise positioning can be achieved via the UWB technology. Since the RSU has the position and movement data of the pedestrian and of the nearby traffic, the received CPM and CAM messages allow the algorithm (running on the smartphone of the pedestrian) to compute possible collisions as a result of this scenario. In such case, the pedestrian is warned via the (wearable) HMI. A second implemented mitigating action in this scenario is that the RSU requests the traffic light management system to turn the traffic lights to a red signal. With these two mitigating actions in place, a possible collision between a pedestrian and a hijacked vehicle, or even normal vehicles, is avoided.
\subsection{Hacked internal vehicle signals}~\label{sec:scen5}
\begin{figure}
    \color{white}\renewcommand\color[2][]{}
        \def\svgwidth{0.9\columnwidth}
\begingroup%
  \makeatletter%
  \providecommand\color[2][]{%
    \errmessage{(Inkscape) Color is used for the text in Inkscape, but the package 'color.sty' is not loaded}%
    \renewcommand\color[2][]{}%
  }%
  \providecommand\transparent[1]{%
    \errmessage{(Inkscape) Transparency is used (non-zero) for the text in Inkscape, but the package 'transparent.sty' is not loaded}%
    \renewcommand\transparent[1]{}%
  }%
  \providecommand\rotatebox[2]{#2}%
  \newcommand*\fsize{\dimexpr\f@size pt\relax}%
  \newcommand*\lineheight[1]{\fontsize{\fsize}{#1\fsize}\selectfont}%
  \ifx\svgwidth\undefined%
    \setlength{\unitlength}{1239bp}%
    \ifx\svgscale\undefined%
      \relax%
    \else%
      \setlength{\unitlength}{\unitlength * \real{\svgscale}}%
    \fi%
  \else%
    \setlength{\unitlength}{\svgwidth}%
  \fi%
  \global\let\svgwidth\undefined%
  \global\let\svgscale\undefined%
  \makeatother%
  \begin{picture}(1,0.62953995)%
    \lineheight{1}%
    \setlength\tabcolsep{0pt}%
    \put(0,0){\includegraphics[width=\unitlength,page=1]{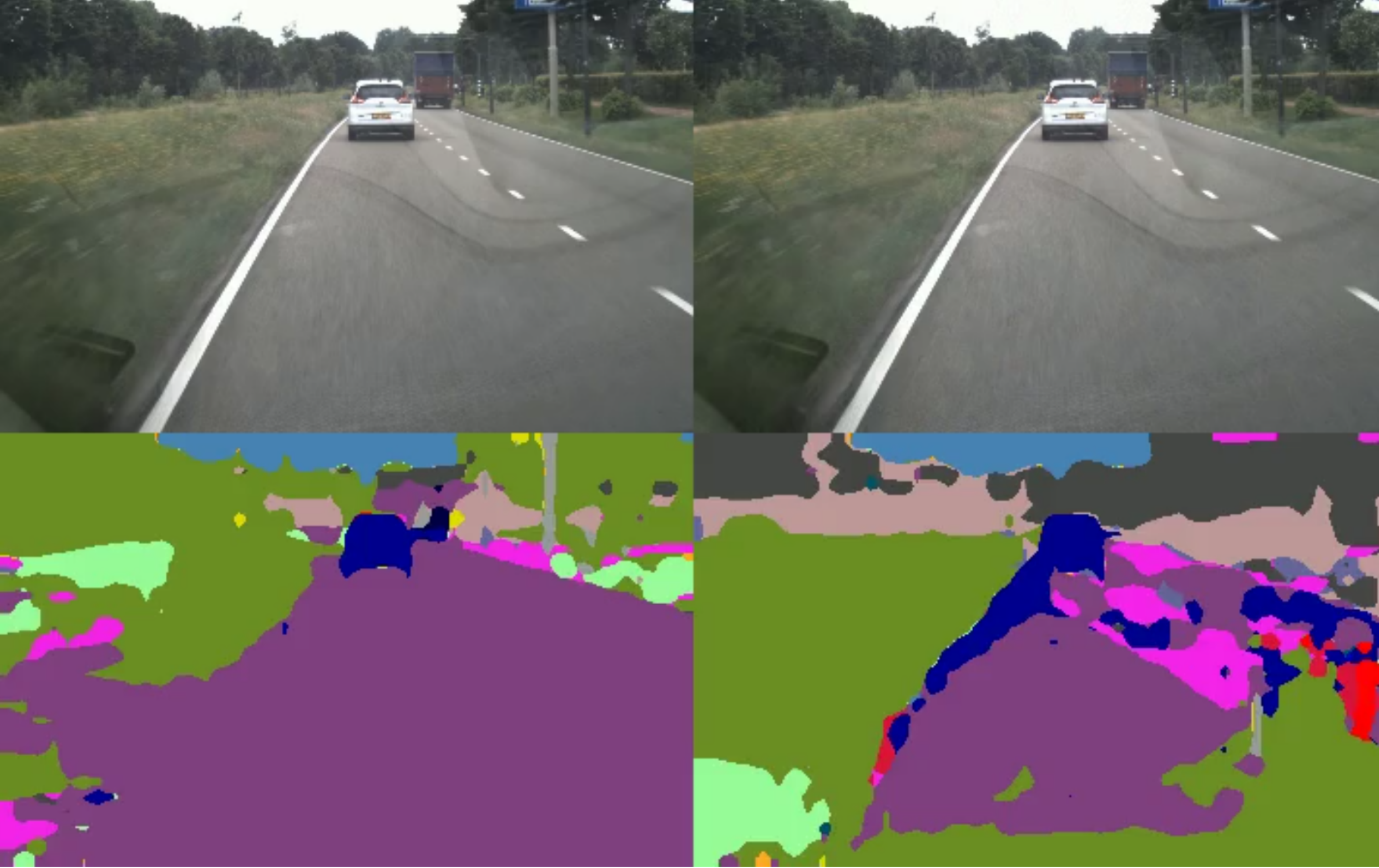}}%
    \put(0.20173536,0.30126805){\color[rgb]{0,0,0}\makebox(0,0)[lt]{\lineheight{1.25}\smash{\begin{tabular}[t]{l}Original feed\end{tabular}}}}%
    \put(0.71188621,0.29968325){\color[rgb]{0,0,0}\makebox(0,0)[lt]{\lineheight{1.25}\smash{\begin{tabular}[t]{l}Hacked feed\end{tabular}}}}%
    \put(0,0){\includegraphics[width=\unitlength,page=2]{fig10.pdf}}%
  \end{picture}%
\endgroup%

    \caption{Original feed and semantic segmentation on the left side versus the feed and segmentation resulting from an adverserial attack on the right.}
    \label{fig:hackedimages}
\end{figure}
On-board cameras are used to detect objects in the environment of the autonomous vehicle. While the autonomous vehicle is driving, several algorithms depend on the sensors of the car. For example, an object detection algorithm identifies items in the video footage coming from the on-board cameras. Image segmentation algorithms based on artificial neural networks deliver a granular understanding of the entire image, as they provide information for every single pixel. Malicious sensor spoofing represents a significant threat to autonomous vehicles. External attackers may try to hack an autonomous driving vehicle approaching an intersection by gaining access to the input-image-stream of a vehicle and targeting the front on-board camera sensor with manipulated images to confuse the algorithm that identifies objects in the video footage, aiming for the vehicle to crash. Usually, the image is manipulated only slightly with the aim that the attack cannot be detected but the effect of the manipulation on the algorithms might be profound~(as demonstrated in Fig.~\ref{fig:hackedimages}).

In order to avoid this situation, a novel unifying method of detecting both anomalous and adversarial inputs by observing activations from the image segmentation model’s layers under attack and transforming these activations using normalising flows onto a target distribution of choice has been developed. It works by first measuring the activations, similar to an EEG on humans, and secondly transforming these observations and lastly using a simple classifier to generate the output - is it anomalous or not? In comparison to many current methods, measurements takes place on all layers, a full brain scan - so to say, not just observing the final outputs. This not only makes it more agnostic to the task at hand but also arguably more robust against potential attackers, as not only the segmentation model has to be fooled but also the detection model at multiple instances~\cite{DAAIN}.

Successfully identified attacks by the anomaly detector are signaled to the on-board unit (OBU) by sending a notification message. Once the OBU receives the signal that something is wrong, it starts sending DENM messages via V2X communication, indicating that  something is wrong. Subsequently the vehicle flushes its certificates. The equipment on the side of the road (such as cameras or RSUs) as well as vehicles around the hacked vehicle, understand that they can not rely on communication from the hacked vehicle because this vehicle itself has concluded that something is wrong.

A device-independent methodology for detecting anomalies  inside the vehicle could be related to checking the integrity of the inter-device communication. In vehicles, multiple communication technologies are used. The most widespread is the Controller Area Network (CAN), which allows the embedded Electronic Control Units (ECUs) to communicate. CAN is a message-based bus communication protocol enabling every participating device to send or receive any message. The CAN network was designed to handle safety problems but not security ones. Every message contains a checksum to detect randomly occurring errors during transmission. However, it is not a message authentication code, thus it does not protect against an intelligent attack.  The protocol neither supports a sender nor a receiver address, and there is no authentication. Messages only contain an identifier that signals the type of content.

Many researchers have previously pointed out these potential security issues, and it was shown in practice that attacks are indeed possible. There are multiple ways an attacker can interfere with the CAN bus: 
\begin{itemize}
    \item The ID fields also serve as a priority identifier. The bus can be flooded with messages of the highest priority preventing the transmission of benign traffic.
    \item A specific functionality can be altered by injecting additional messages of that ID into the traffic. These new messages are accepted, as there is no mechanism on the bus the differentiate between the malicious and the benign messages. The number of injected messages needs to be sufficiently high to suppress the original messages because those are also still transmitted.
    \item Upon a compromise of an ECU, an attacker could alter the contents of an original message. This approach does not result in additional messages appearing on the bus; however, achieving such compromise is more complicated.
\end{itemize}

Many potential detection approaches have been published in the literature. The mostly regular repetition times of CAN messages make it relatively simple to detect message injection attacks. A small portion of the potential ID values is used only, making it trivial to detect a DoS attack containing the highest priority messages.
Other message injection attacks can be detected based on the changes in the regular repetition times~\cite{DBLP:conf/esorics/GazdagNBS18}. Practical solutions have been developed with both statistical and machine learning based approaches.
The most complex detection challenge is the message data modification attack~\cite{DBLP:conf/vehits/ChiscopG0B21}. There is no additional meta-information coming from the transmission in this scenario, so the detection has to rely solely on the received data. There are effective anomaly detection solutions proposed for this scenario as well; however, none of them are as efficient as the detectors against message injection attacks.

\section{Conclusion and Future Work}
Communication with other traffic participants could bring high safety benefits to automated driving. However, these communication channels have certain inherent vulnerabilities, which could impose great risks. In this demonstrator, we have made use of external sensor data from digital infrastructure for anomaly and misbehavior detection. Based on cooperative awareness services, we generate warnings using DENM messages (which are available at a central level) as mitigating actions to cybersecurity threats. We succesfully use traffic light changes as a mitigating action on detected anomalies or misbehaviors and, finally, we have extended current CAM-only misbehavior reporting ~\cite{ETSI103759} with CPM-based anomaly detection and DENM warning messaging. Not only is this principle used for flagging and mitigating cybersecurity threats, the same principle is used for an internal attack (as demonstrated in section 3.5). In the SECREDAS project, we have succesfully realized proof-of-concept demonstrators for these techniques using a cooperative intersection and cooperative vehicles~\footnote{The reader is deferred to \hyperlink{https://secredas-project.eu/}{\tt https://secredas-project.eu/} for full coverage of the project scope. The mentioned video material is shown at~ \hyperlink{https://player.secredas-project.eu/}{\tt https://player.secredas-project.eu/}.}. \\\\
As future work, next to the performed functional demonstrations, we would like to extend our analysis with performance assessment of the concepts in a quantitative manner, i.e., through definition of key performance indicators, to measure the wide impact on safety of our methodologies. Other interesting future work includes extending for more relevant scenarios. In the scope of this work, only attacks on a cooperative intersection have been considered. Although we have not covered all the potential attacks in such an environment, it would be interesting to look at the new challenges that appear in other relevant environments. Moreover, we see opportunities to extend the number of mitigating actions and improve on them. Namely, the mitigating actions need to be deployed in such a way that they improve the safety and security of CCAM, while not hampering the benefits of connected driving (which could lead to a loss of human acceptance for the technology, hence drastically reducing the commercial viability). For example, handling an attack detected on the CAN level is unclear. This communication bus with most of its participating ECUs is essential to the operation of the vehicle, therefore shutting down the bus or blocking the communication of specific ECUs is not an option. A potential reaction could be to switch to a reduced operation, which for example, limits the vehicle's speed. Another possible follow-up step is to stop the V2X communication as the transmitted data is potentially corrupted.

\section{Acknowledgements}
All partners have contributed to the writing of the paper. TNO Integrated Vehicle Safety hosted and coordinated the demonstration activities and performed the research and implementation of the anomaly detection algorithm of Section~\ref{sec:scen1}. Canon Research Centre France has performed the research on the anomaly detection algorithms of Section~\ref{sec:scen2}. Commsignia peformed the research on the traffic light anomaly detection algorithm of Section~\ref{sec:scen3} and the UWB-based collision detection algorithm of Section~\ref{sec:scen4}. Merantix Momentum GmbH has performed the research on the anomaly detection for onboard camera attacks in Section~\ref{sec:scen5}. CrySyS Lab has performed the research on the CAN anomaly detection algorithms from Section~\ref{sec:scen5}\\\\ 
This work has been funded by EU ECSEL Project SECREDAS Cyber Security for Cross Domain Reliable Dependable Automated Systems (Grant Number 783119). The authors would like to thank the consortium for the successful cooperation.

\errorcontextlines=99

\end{document}